\documentclass[reqno,centertags,a4paper,10pt]{amsart}

\usepackage{amssymb,amsmath,amsthm}

\newcommand{\R}{\mathbb{R}} 
\newcommand{\T}{\mathbb{T}} \newcommand{\C}{\mathbb{C}}
\newcommand{\N}{\mathbb{N}}

\theoremstyle{plain}

\newtheorem{lemma}{Lemma}
\newtheorem{kor}{Corollary}
\newtheorem{satz}{Theorem}

\newcommand{\F}{\mbox{${\mathcal F}$}}
\newcommand{\eps}{\varepsilon}
\newcommand{\lb}{\langle} \newcommand{\rb}{\rangle}
\newcommand{\ls}{\lesssim}\newcommand{\gs}{\gtrsim}

\begin{document}

\thanks{M.P. is partially supported by the Funda\c{c}\~ao para a Ci\^encia e a 
 Tecnologia through the program POCI 2010/FEDER and by the grant SFRH/BPD/22018/2005}
\thanks{J.D.S. is partially supported by the Funda\c{c}\~ao para a Ci\^encia e a 
 Tecnologia through the program POCI 2010/FEDER and by the project POCI/FEDER/MAT/55745/2004}
\author[A.~Gr{\"u}nrock]{Axel~Gr{\"u}nrock} 
\author[M.~Panthee]{Mahendra~Panthee}
\author[J.D.~Silva]{Jorge~Drumond~Silva}
\title[GWP below $L^2$ for gKdV-3]
{A remark on global well-posedness below $L^2$ for the gKdV-3 equation}

\address{Axel~Gr{\"u}nrock: Bergische Universit\"at
  Wuppertal, Fachbereich C: Mathematik / Naturwissenschaften,
  Gau{\ss}stra{\ss}e 20, 42097 Wuppertal, Germany.}
\email{axel.gruenrock@math.uni-wuppertal.de}
\address{Mahendra~Panthee:
 Centro de An\'alise Matem\'atica, Geometria e Sistemas Din\^amicos, 
 Departamento de Matem\'atica, Instituto Superior T\'ecnico, 
 1049-001 Lisboa, Portugal.}
\email{mpanthee@math.ist.utl.pt}
\address{Jorge~Drumond~Silva:
 Centro de An\'alise Matem\'atica, Geometria e Sistemas Din\^amicos, 
 Departamento de Matem\'atica, Instituto Superior T\'ecnico, 
 1049-001 Lisboa, Portugal.}
\email{jsilva@math.ist.utl.pt}

\begin{abstract}
The $I$-method in its first version as developed by Colliander et al. in 
\cite{Q1} is applied to prove that the Cauchy-problem for the generalised 
Korteweg-de Vries equation of order three (gKdV-3) is globally well-posed 
for large real-valued data in the Sobolev space $H^s(\R \rightarrow \R)$, provided 
$s>-\frac{1}{42}$.
\end{abstract}

\maketitle
\section{Introduction}

In a recently published paper of Tao \cite{T07} concerning the Cauchy-problem
for the generalised Korteweg-de Vries equation of order three (for short: gKdV-3), 
i.e.:
\begin{equation}\label{gkdv3}
u_t + u_{xxx} \pm (u^4)_x=0,\hspace{1cm}u(0,x)=u_0(x) , \,\,\, x \in {\R},
\end{equation}
it was shown that this problem is locally well-posed for data $u_0$ in the 
critical Sobolev space $\dot{H}^{-\frac{1}{6}}(\R \rightarrow \C)$ and globally 
well-posed for data with sufficiently small $\dot{H}^{-\frac{1}{6}}$-norm. 
Moreover, scattering results in $H^1 \cap \dot{H}^{-\frac{1}{6}}(\R \rightarrow \R)$ 
for the radiation component of a perturbed soliton were obtained. Tao's local 
result improves earlier work of Kenig, Ponce, and Vega ($s\ge \frac{1}{12}$, 
see \cite[Theorem 2.6]{KPV93}) and of the first author ($s>-\frac{1}{6}$, see \cite{G05}), 
while the global \emph{small} data theory seems to be completely new in Sobolev spaces 
of negative index. For \emph{large real valued} data  $u_0 \in H^s(\R \rightarrow \R)$, 
$s \ge 0$, global well-posedness of \eqref{gkdv3} was obtained in \cite{G05} by 
combining the conservation of the $L^2$-norm with the local $L^2$-result, for 
$s\ge 1$ this was already in \cite[Corollary 2.7]{KPV93}, where the energy 
conservation was used. 

\vspace{0.3cm}

Starting with Bourgain's splitting argument \cite{B98} and followed by 
the ``$I$-method'' or ``method of almost conservation laws''
introduced and further refined by Colliander, Keel, Staffilani, Takaoka, and 
Tao in a series of papers - see e. g. \cite{Q1}, \cite{Q2}, \cite{Q3}, \cite{Q4},
 \cite{Q5} - effective techniques have been developed, which are capable to show 
large data global well-posedness \emph{below} certain conserved quantities such 
as the energy or the $L^2$-norm. The question of whether and to what extent these 
methods apply to the Cauchy-problem for gKdV-3, was raised as well by Linares 
and Ponce \cite[p.177 and p.183]{LP04} as by Tao, see Remark 5.3 in \cite{T07}. 
In this note we establish global well-posedness of \eqref{gkdv3} for large data 
$u_0 \in H^s(\R \rightarrow \R)$, provided $s > - \frac{1}{42}$, thus giving a 
partial answer to this question. Our proof combines the first version of the 
$I$-method as in \cite{Q1} with a sharp four-linear $X_{s,b}$-estimate exhibiting 
an extra\footnote{i. e. beyond the cancellation of the derivative in the 
nonlinearity} gain of half a derivative. 

\vspace{0.3cm}

Before we turn to the details, let us point out that 
substantial difficulties appear, if we try to push the analysis further to 
lower values of $s$; by following the construction of a sequence of ``modified 
energies'' in \cite{Q4} we are led already in the second step to a Fourier 
multiplier, say $\mu_8$, corresponding to $M_4$ in \cite{Q4}, with a quadratic 
singularity, and the argument breaks down.\footnote{A similar problem was observed 
by Tzirakis for the quintic semilinear Schr\"odinger equation in one dimension, 
see the concluding remark in \cite{Tz05}.} Our fruitless effort in this direction 
seems to confirm Tao's remark, that ``it is unlikely that these methods would 
get arbitrarily close to the scaling regularity $s=-\frac{1}{6}$.'' \cite[Remark 5.3]{T07}

\vspace{0.3cm}

{\bf{Acknowledgement:}} The first author, A. G., wishes to thank the Center of 
Mathematical Analysis, Geometry and Dynamical Systems at the IST in Lisbon for 
its kind hospitality during his visit.

\section{A variant of local well-posedness, the decay estimate, and the main result}
Here we follow the lines of \cite{Q1}: The operator $I_N$ is defined via the Fourier transform by
\[\widehat{I_Nu} (\xi):= m(\frac{|\xi|}{N})\widehat{u}(\xi),\]
where $m:\R^+ \rightarrow \R^+$ is a smooth monotonic function with $m(x)=1$ for 
$x\le 1$ and $m(x)=x^s$, $x\ge 2$. Here $s<0$, so that $0<m(x)\le 1$. 
$I_N : H^s \rightarrow L^2$ is isomorphic and $\|I_N \cdot\|_{L^2}$ 
defines an equivalent norm on $H^s$, with implicit constants depending on $N$. 

\vspace{0.3cm}

The crucial nonlinear estimate in the proof of local well-posedness for 
\eqref{gkdv3} with $H^s$-data, $s>-\frac{1}{6}$, is
\begin{equation}\label{gruen}
\|\partial_x \prod_{i=1}^4 u_i\|_{X_{s,b'}} \ls \prod_{i=1}^4 \|u_i\|_{X_{s,b}},
\end{equation}
which holds true, whenever $0\ge s > -\frac{1}{6}$, $-\frac{1}{2}<b'<s-\frac{1}{3}$
and $b>\frac{1}{2}$, see \cite[Theorem 1]{G05}. The $X_{s,b}$-norms used here are 
given by
\[\|u\|_{X_{s,b}}=\left(\int d \xi d\tau \lb\tau - \xi^3\rb^{2b}
 \lb\xi\rb^{2s}|\F u(\xi,\tau)|^2 \right) ^{\frac{1}{2}},\]
where $\F$ denotes the Fourier transform in both variables. Later on we shall also 
use the restriction norms 
$\|v\|_{X_{s,b}(\delta)}=\inf{\{\|u\|_{X_{s,b}}: u|_{[0,\delta]\times \R}=v\}}$.
Applying the interpolation lemma \cite[Lemma 12.1]{Q5} to \eqref{gruen} we obtain,
under the same assumptions on the parameters $s$, $b'$ and $b$,
\begin{equation}\label{vargruen}
\|I_N\partial_x \prod_{i=1}^4 u_i\|_{X_{0,b'}} \ls \prod_{i=1}^4 \|I_Nu_i\|_{X_{0,b}},
\end{equation}
where the implicit constant is \emph{independent} of $N$. Now familiar arguments 
invoking the contraction mapping principle give the following variant of local 
well-posedness.

\begin{lemma}\label{varlwp}
For $s>-\frac{1}{6}$ the Cauchy-problem \eqref{gkdv3} is locally well-posed 
for data $u_0 \in (H^s,\|I_N\cdot\|_{L^2})$. The lifespan $\delta$ of the local solution 
$u$ satisfies
\begin{equation}\label{1.5}
\delta \gs \|I_Nu_0\|_{L^2}^{-\frac{18}{6s+1}-}
\end{equation}
and moreover we have for $b=\frac{1}{2}+$
\begin{equation}\label{1.6}
\|I_Nu\|_{X_{0,b}(\delta)}\ls \|I_Nu_0\|_{L^2}.
\end{equation}
\end{lemma}

Replacing $u^2$ by $u^4$ in the calculation on p. 2 of \cite{Q1}, 
we obtain for a solution $u$ of \eqref{gkdv3}
\begin{equation}\label{1.3/1.4}
\|I_Nu(\delta)\|^2_{L^2}-\|I_Nu(0)\|^2_{L^2} \ls 
 \|\partial_x (I_Nu^4-(I_Nu)^4)\|_{X_{0,-b}(\delta)}\|I_Nu\|_{X_{0,b}(\delta)}.
\end{equation}
The next section will be devoted to the proof that for $b>\frac{1}{2}$, $0\ge s \ge -\frac{1}{8}$
\begin{equation}\label{decay1}
\|\partial_x (I_Nu^4-(I_Nu)^4)\|_{X_{0,-b}(\delta)}
 \ls N^{-\frac{1}{2}}\|I_Nu(0)\|_{L^2}^4
\end{equation}
(see Corollary \ref{last} below), which together with \eqref{1.3/1.4} and \eqref{1.6} gives
\begin{equation}\label{decay}
\|I_Nu(\delta)\|_{L^2}-\|I_Nu(0)\|_{L^2} \ls N^{-\frac{1}{2}}\|I_Nu(0)\|_{L^2}^4.
\end{equation}
Now the decay estimate \eqref{decay} allows us to prove our main result:
\begin{satz}\label{main}
Let $s>-\frac{1}{42}$ and $u_0 \in H^s(\R \rightarrow \R)$. Then the solution 
$u$ of \eqref{gkdv3} according to Lemma \ref{varlwp} extends uniquely to any time interval 
$[0,T]$ and satisfies
\begin{equation}\label{gb}
\sup_{0 \le t \le T}\|u(t)\|_{H^s} \ls \lb T \rb^{\frac{-2s}{1+42s}}\|u_0\|_{H^s}.
\end{equation}
\end{satz}

Proof: We choose $\eps_0$ so that Lemma \ref{varlwp} gives the lifespan $\delta =1$ 
for all data $\phi \in H^s$ with $\|I_N\phi\|_{L^2} \le 2 \eps_0$. Moreover we 
demand $16 C \eps_0^3 \le 1$, where $C$ is the implicit constant in the decay 
estimate \eqref{decay}. Assuming without loss that $T \gg 1$, we fix parameters 
$C_1$, $N$ and $\lambda$ with
\[2C_1^{-\frac{1}{6}-s}\|u_0\|_{H^s}=\eps_0,\hspace{0.5cm}
 N^{\frac{1+42s}{2(1+6s)}}=C_1^3 T,\hspace{0.5cm}\mbox{and}\hspace{0.5cm}
 \lambda=C_1N^{\frac{-6s}{1+6s}}.\]
Then $N^{\frac{1}{2}}=\lambda^3T$ and for $u_0^{\lambda}(x)=\lambda^{-\frac{2}{3}}u_0(\frac{x}{\lambda})$ 
it is easily checked that $\|I_Nu_0^{\lambda}\|_{L^2}\le \eps_0$. 
For any $k \in \N$, repeated applications of Lemma \ref{varlwp} give a solution $u^{\lambda}$ of gKdV-3 
with $u^{\lambda}(0)=u_0^{\lambda}$ on $[0,k]$, as long as
\begin{equation}\label{cond}
\|I_Nu^{\lambda}(j)\|_{L^2}\le 2 \eps_0, \qquad 1 \le j < k.
\end{equation}
Since by \eqref{decay} and the second assumption on $\eps_0$
\[\|I_Nu^{\lambda}(j)\|_{L^2}\le \eps_0 + jCN^{-\frac{1}{2}}(2\eps_0)^4
 \le \eps_0(1 + jN^{-\frac{1}{2}}),\]
condition \eqref{cond} is fulfilled for $j \le N^{\frac{1}{2}}=\lambda^3T$. Thus 
$u^{\lambda}$ is defined on $[0,\lambda^3T]$, and with 
$u(x,t)=\lambda^{\frac{2}{3}}u^{\lambda}(\lambda x,\lambda^3t)$ we obtain a 
solution of \eqref{gkdv3} on $[0,T]$. Finally we have for $0 \le t \le T$
\[\|u(t)\|_{H^s} \ls \|I_{\lambda N}u(t)\|_{L^2} \ls 
 \lambda^{\frac{1}{6}}\|I_N u^{\lambda}(\lambda^3t)\|_{L^2}\]
with $\lambda^{\frac{1}{6}} \sim T^{\frac{-2s}{1+42s}}$ and 
$\|I_N u^{\lambda}(\lambda^3t)\|_{L^2}$ being bounded during the whole 
iteration process by $2\eps_0 \ls \|u_0\|_{H^s}$. This gives the growth bound \eqref{gb}.
$\hfill \Box$

\section{The decisive four-linear estimate}

Let us first recall several linear and bilinear Airy estimates (in their 
$X_{s,b}$-versions), which shall be used below; by interpolation between 
the sharp version of Kato's smoothing effect (see \cite[Theorem 4.1]{KPV91}) 
and the maximal function estimate from \cite[Theorem 3]{S87} we have
\begin{equation}\label{lin1}
\|J^s u\|_{L^p_x(L^q_t)} \ls \|u\|_{X_{0,b}}, 
\end{equation}
whenever $b>\frac{1}{2}$, $-\frac{1}{4} \le s \le 1$ and 
$(\frac{1}{p},\frac{1}{q}) = (\frac{1-s}{5},\frac{1+4s}{10})$. We will use 
\eqref{lin1} with $s=0$, i.e.
\begin{equation}\label{lin2}
\| u\|_{L^5_x(L^{10}_t)} \ls \|u\|_{X_{0,b}}, 
\end{equation}
and the dual version of \eqref{lin1} with $s=\frac{1}{2}$, which is
\begin{equation}\label{lin3}
\| u\|_{X_{\frac{1}{2},-b}} \ls \|u\|_{L^{\frac{10}{9}}_x(L^{\frac{10}{7}}_t)}. 
\end{equation}
Moreover we shall rely on the Strichartz type estimate
\begin{equation}\label{lin4}
\| u\|_{L^8_{xt}} \ls \|u\|_{X_{0,b}}, \,\,\,\,(b>\frac{1}{2})
\end{equation}
(cf. \cite[Theorem 2.4]{KPV91}) and the bilinear estimate
\begin{equation}\label{bil}
\| I^{\frac{1}{2}}I_-^{\frac{1}{2}}(u,v)\|_{L^2_{xt}} \ls \|u\|_{X_{0,b}}\|v\|_{X_{0,b}}, \,\,\,\,(b>\frac{1}{2})
\end{equation}
from \cite[Corollary 1]{G05}. Here $I^s$ ($J^s$) denotes the Riesz (Bessel) potential operator of order
$-s$ and $I^s_-$ is defined via the Fourier transform by
\[\widehat{ I_-^s (f,g)} (\xi) := \int_{\xi_1+\xi_2=\xi}d\xi_1|\xi_1-\xi_2|^s \widehat{f}(\xi_1)\widehat{g}(\xi_2).\]
Now we turn to the crucial four-linear $X_{s,b}$-estimate:
\begin{lemma}\label{lemma}
Let $b>\frac{1}{2}$, $s_i \le 0$, $1\le i \le 4$, with $\sum_{i=1}^4 s_i = -\frac{1}{2}$. 
Then
\begin{equation}\label{qual}
\|\partial_x \prod_{i=1}^4 v_i\|_{X_{0,-b}} \ls \prod_{i=1}^4 \|v_i\|_{X_{s_i,b}}.
\end{equation}
\end{lemma}

Proof: We write
\[\|\partial_x \prod_{i=1}^4 v_i\|_{X_{0,-b}} = c
 \|\xi\lb\tau-\xi^3\rb^{-b} \int d\nu \prod_{i=1}^4\F{v_i}(\xi_i,\tau_i)\|_{L^2_{\xi,\tau}},\]
where $d \nu = d\xi_1..d\xi_{3} d\tau_1.. d \tau_{3}$ and $\sum_{i=1}^4 (\xi_i,\tau_i) = (\xi, \tau)$, 
and divide the domain of integration into three regions $A$, $B$ and $C=(A\cup B)^c$.
In region $A$ we assume that\footnote{Here $\xi_{max}$ is defined by 
$|\xi_{max}| = \max_{i=1}^4 |\xi_i|$, similarly $\xi_{min}$.} $|\xi_{max}| \le 1$
and hence $|\xi|\le 4$, so for this region we get the upper bound
\[\|\prod_{i=1}^4 J^{s_i}v_i\|_{L^2_{xt}} \le \prod_{i=1}^4\|J^{s_i}v_i\|_{L^8_{xt}}
  \ls \prod_{i=1}^4\|v_i\|_{X_{s_i,b}},\]
where in the last step we have used the $L^8_{xt}$-Strichartz-type estimate \eqref{lin4}.
Concerning the region $B$ we shall assume - besides $|\xi_{max}| \ge 1$, implying $\lb\xi_{max}\rb \ls |\xi_{max}|$ - that
\begin{itemize}
\item[i)] $|\xi_{min}| \le 0.99 |\xi_{max}|$ or
\item[ii)] $|\xi_{min}| > 0.99 |\xi_{max}|$, and there are exactly two indices 
           $i \in \{1,2,3,4 \}$ with $\xi_i > 0$.
\end{itemize}
Then the region $B$ can be split further into a finite number of subregions, so that 
for any of these subregions there exists a permutation $\pi$ of $\{1,2,3,4 \}$ with
\[|\xi| \ls |\xi|^{\frac{1}{2}}|\xi_{\pi(1)}+\xi_{\pi(2)}|^{\frac{1}{2}}
 |\xi_{\pi(1)}-\xi_{\pi(2)}|^{\frac{1}{2}}\prod_{i=1}^4\lb\xi_i\rb^{s_i}.\]
Assume $\pi = id$ for the sake of simplicity now. Then we get the upper bound
\begin{eqnarray*}
&& \|(I^{\frac{1}{2}}I_-^{\frac{1}{2}}(J^{s_1}v_1,J^{s_2}v_2))(J^{s_3}v_3)(J^{s_4}v_4)\|_{X_{\frac{1}{2},-b}} \\
&\ls & \|(I^{\frac{1}{2}}I_-^{\frac{1}{2}}(J^{s_1}v_1,J^{s_2}v_2))(J^{s_3}v_3)(J^{s_4}v_4)\|_{L^{\frac{10}{9}}_x(L^{\frac{10}{7}}_t)}  \\
&\ls &  \|I^{\frac{1}{2}}I_-^{\frac{1}{2}}(J^{s_1}v_1,J^{s_2}v_2)\|_{L^2_{xt}}\|J^{s_3}v_3\|_{L^5_x(L^{10}_t)}\|J^{s_4}v_4\|_{L^5_x(L^{10}_t)} 
\ls \prod_{i=1}^4\|v_i\|_{X_{s_i,b}}.
\end{eqnarray*}
Here we have applied the estimates \eqref{lin3}, H\"older, \eqref{bil} and \eqref{lin2}. 
Finally we consider the remaining region $C$: Here the $|\xi_i|$, $1\le i \le 4$, are 
all very close together and $\gs  \lb\xi_i\rb$. Moreover, at least three of the variables 
$\xi_i$ have the same sign. Thus for the quantity $c.q.$ controlled by the expressions 
$\lb\tau - \xi ^3\rb$, $\lb\tau_i - \xi_i ^3\rb$, $1\le i \le 4$, we have in this region:
\[c.q.:=|\xi ^3 - \sum_{i=1}^4 \xi_i ^3| \gs \sum_{i=1}^4 \lb\xi_i \rb^3 \gs \lb\xi \rb^3.\]
So the contribution of the subregion, where $\lb\tau - \xi ^3\rb \ge \max_{i=1}^4{\lb\tau_i - \xi_i ^3\rb}$,
is bounded by
\[\|\prod_{i=1}^4 J^{s_i}v_i\|_{L^2_{xt}} \le \prod_{i=1}^4\|J^{s_i}v_i\|_{L^8_{xt}}
  \ls \prod_{i=1}^4\|v_i\|_{X_{s_i,b}},\]
where \eqref{lin4} was used again. On the other hand, if $\lb\tau_1 - \xi_1 ^3\rb$ is 
dominant, we write $\Lambda^{\frac{1}{2}}=\F^{-1}\lb\tau - \xi ^3\rb^{\frac{1}{2}}\F$ 
and obtain the upper bound
\begin{eqnarray*}
 \|(\Lambda^{\frac{1}{2}}J^{s_1}v_1)\prod_{i=2}^4J^{s_i}v_i\|_{X_{0,-b}} 
&\ls & \|(\Lambda^{\frac{1}{2}}J^{s_1}v_1)\prod_{i=2}^4J^{s_i}v_i\|_{L^{\frac{8}{7}}_{xt}}\\
\le  \|\Lambda^{\frac{1}{2}}J^{s_1}v_1\|_{L^2_{xt}}\prod_{i=2}^4\|J^{s_i}v_i\|_{L^8_{xt}}
& \ls & \prod_{i=1}^4\|v_i\|_{X_{s_i,b}}.
\end{eqnarray*}
Here the dual version $X_{0,-b} \supset L^{\frac{8}{7}}_{xt}$ of the $L^8_{xt}$ estimate was 
used first, followed by H\"older's inequality and the estimate itself. The remaining subregions, 
where $\lb\tau_k - \xi_k ^3\rb$, $2 \le k \le 4$, are maximal, can be treated in precisely the 
same manner.

$\hfill \Box$

\begin{kor}\label{last}
Let $b> \frac{1}{2}$ and $0\ge s \ge -\frac{1}{8}$. Then
\begin{equation}\label{qualN}
\|\partial_x (I_N (\prod_{i=1}^4u_i)-\prod_{i=1}^4I_Nu_i)\|_{X_{0,-b}(\delta)}
 \ls N^{-\frac{1}{2}}\prod_{i=1}^4\|I_Nu_i\|_{X_{0,b}(\delta)}.
\end{equation}
Especially, if $u$ is a solution of gKdV-3 according to Lemma \ref{varlwp} with 
$u(0)=u_0$, then
\begin{equation}\label{quarticN}
\|\partial_x (I_Nu^4-(I_Nu)^4)\|_{X_{0,-b}(\delta)}
 \ls N^{-\frac{1}{2}}\|I_Nu_0\|_{L^2}^4.
\end{equation}
\end{kor}

Proof: By \eqref{1.6} the estimate \eqref{qualN} implies \eqref{quarticN}. Thus it suffices 
to show
\begin{equation}\label{qualN-d}
\|\partial_x (I_N (\prod_{i=1}^4u_i)-\prod_{i=1}^4I_Nu_i)\|_{X_{0,-b}}
 \ls N^{-\frac{1}{2}}\prod_{i=1}^4\|I_Nu_i\|_{X_{0,b}}.
\end{equation}
Now let $\xi_i$ denote the frequencies of the $u_i$, $1\le i \le 4$. If all the $|\xi_i|\le N$, 
then either $|\xi|\le N$, such that there's no contribution at all, or we have $|\xi|\ge N$, so 
that at least, say, $|\xi_1|\ge \frac{N}{4}$. In this case, by Lemma \ref{lemma}, the norm 
on the left of \eqref{qualN-d} is bounded by
\[\|\partial_x \prod_{i=1}^4 u_i\|_{X_{0,-b}} \ls \|u_1\|_{X_{-\frac{1}{2},b}}
 \prod_{i=2}^4 \|u_i\|_{X_{0,b}} \ls N^{-\frac{1}{2}}\prod_{i=1}^4\|I_Nu_i\|_{X_{0,b}}.\]
Otherwise there are $k$ large frequencies for some $1\le k \le 4$. By symmetry we may 
assume that $|\xi_1|,\dots,|\xi_k|\ge N$ and $|\xi_{k+1}|,\dots,|\xi_4|\le N$. Then we 
have, again by Lemma \ref{lemma},
\begin{eqnarray*}
\|\partial_x \prod_{i=1}^4 I_Nu_i\|_{X_{0,-b}} & \ls & 
 \prod_{i=1}^k \|I_Nu_i\|_{X_{-\frac{1}{2k},b}}\prod_{i=k+1}^4\|I_Nu_i\|_{X_{0,b}} \\
 & \ls & N^{-\frac{1}{2}}\prod_{i=1}^4\|I_Nu_i\|_{X_{0,b}}
\end{eqnarray*}
as well as
\[\|\partial_x I_N\prod_{i=1}^4 u_i\|_{X_{0,-b}} \ls 
 \prod_{i=1}^k \|u_i\|_{X_{-\frac{1}{2k},b}}\prod_{i=k+1}^4\|u_i\|_{X_{0,b}}.\]
Since for any $s_1 \le s$ and any $v$ with frequency $|\xi| \ge N$ it holds that
\[\|v\|_{X_{s_1,b}} \ls N^{s_1-s}\|v\|_{X_{s,b}}\sim N^{s_1}\|I_Nv\|_{X_{0,b}},\]
the latter is again bounded by the right hand side of \eqref{qualN-d}.
$\hfill \Box$

\end{document}